\documentclass[12pt]{article}
\usepackage{fancyhdr}
\usepackage{fullpage}
\usepackage{amsmath}
\usepackage{amssymb}
\usepackage{theorem}
\newtheorem{thm}{Theorem}[section]
\newtheorem{prop}{Proposition}[section]
\newtheorem{corr}{Corollary}[section]
\newtheorem{lemma}{Lemma}[section]
{\theorembodyfont{\rm} 
\newtheorem{rem}{Remark}[section]}
\newtheorem{defe}{Definition}[section]

\numberwithin{equation}{section}
\def\n{{\mathfrak n}}

\def\sp{\text{Span}}
\def\qs{$\quad \square$}
\def\C{\mathbb C}
\def\l{{\mathfrak{l}}}
\def\h{\mathfrak{h}}
\def\g{\mathfrak{g}}

\def\o{\otimes}
\def\Ad{\text{\rm Ad}\:}

\def\End{\text{End}}

\def\R{{\cal R}}
\def\J{{\cal J}}
\def\N{{\mathbb N}}

\begin{document}

\title{Explicit quantization of dynamical r-matrices 
for finite dimensional semisimple Lie algebras}

\author{Pavel Etingof, Travis Schedler, and Olivier Schiffmann}
\date{}
\maketitle

\section{Introduction} \label{is}

\subsection{Classical r-matrices}
In the early eighties, Belavin and Drinfeld \cite{BD}
classified nonskewsymmetric classical \linebreak r-matrices
for simple Lie algebras. It turned out that such r-matrices,  
up to isomorphism and twisting by elements from 
the exterior square of the Cartan subalgebra, are classified by 
combinatorial objects which are now 
called Belavin-Drinfeld triples. By definition, a Belavin-Drinfeld triple 
for a simple Lie algebra $\g$ is 
a triple $(\Gamma_1,\Gamma_2,T)$, where $\Gamma_1,\Gamma_2$ are subsets 
of the Dynkin diagram $\Gamma$ of $\g$, and $T:\Gamma_1\to \Gamma_2$ is 
an isomorphism which preserves the inner product and satisfies the nilpotency 
condition: if $\alpha\in \Gamma_1$ then there exists $k$ such that 
$T^{k-1}(\alpha)\in \Gamma_1$ but $T^k(\alpha)\notin \Gamma_1$. 
The r-matrix corresponding to such a triple is given by 
a certain explicit formula. These results generalize in a
straightforward way to semisimple Lie algebras. 

In \cite{S}, the third author generalized the work of 
Belavin and Drinfeld and 
classified  classical
nonskewsymmetric \textit{dynamical} r-matrices for simple Lie algebras. 
It turns out that they have an even simpler classification: up 
to gauge transformations, they are classified by generalized 
Belavin-Drinfeld triples, which are defined as the usual Belavin-Drinfeld 
triples but without any nilpotency condition.  
The dynamical r-matrix corresponding to such a triple is given by 
a certain explicit formula. As before, these results can be
generalized to semisimple Lie algebras.

\subsection{Quantization of r-matrices} The problem of quantization of the Belavin-Drinfeld r-matrices  
(i.e. finding the corresponding quantum R-matrices) 
has been open for a long time. The history of this problem is as
follows. 

In the case when $\Gamma_i$ are empty
(the ``standard'' r-matrix), the quantization was
provided by Drinfeld and Jimbo in mid 80-s, which gave birth 
to the modern theory of quantum groups.
 
In 1990, Cremmer and Gervais \cite{CG} provided an explicit
quantum R-matrix
(in the vector representation) for the so-called Cremmer-Gervais
triple 
for $sl_n$ -- the Belavin-Drinfeld
triple where $\Gamma_1$ and $\Gamma_2$ are the whole
Dynkin diagram without the last and the first vertex,
respectively, 
and $T$ is the shift to the right by one position. 
Other proofs of the fact that the Cremmer-Gervais R-matrix
satisfies the quantum Yang-Baxter equation 
and the Hecke relation were given 
by Hodges \cite{H1,H2}.

In 1992, Gerstenhaber, Giaquinto,
and Schack \cite{GGS} suggested a conjectural explicit quantization 
of all Belavin-Drinfeld r-matrices for the Lie algebra $sl_n$, in the
vector representation (the GGS conjecture); it has been checked in many cases
\cite{GH,Sch1}, but a general proof is still unavailable. 

In 1995, it was shown in \cite{EK} that a quantization exists in
principle, but the method of \cite{EK} does not allow one to compute
the quantization explicitly.

In 1996, Hodges \cite{H3} suggested an explicit
quantization in the case when $\Gamma_1\cap\Gamma_2=\emptyset$
(for any Lie algebra), which yields a proof of the GGS conjecture in this
case. Namely, he constructed a twist which transforms the 
usual quantum group $U_q(\g)$ to a quantum group that is a
quantization of the Lie bialgebra corresponding to the 
given Belavin-Drinfeld triple with disjoint $\Gamma_1$ and
$\Gamma_2$. 

In early 1999, the second author generalized the method of Hodges 
to a wider class of triples (\cite{Sch2}). He also did computer calculations
which led him to a conjectural quantization 
of all triples for $sl_n$ in the vector representation.  
This work led to an understanding what the generalization 
of Hodges' formula to any
triple should look like, and eventually to a complete 
explicit solution of the problem, which is given here.

\subsection{Description of the 
paper} In this paper, we provide an explicit quantization 
of dynamical r-matrices for semisimple Lie algebras,
classified in \cite{S}, which includes the Belavin-Drinfeld 
r-matrices. We do so by constructing an appropriate (dynamical) twist
in the tensor square of the Drinfeld-Jimbo quantum group $U_q(\g)$. 
The construction of this twist is based on the method 
stemming from \cite{JKOS} and \cite{ABRR}, i.e. on 
defining the twist as a unique solution
of a suitable difference equation. 
This yields 
a simple closed formula for the twist.

In the 
case of ordinary Belavin-Drinfeld triples (i.e. satisfying the 
nilpotency condition), the constructed dynamical twist 
can be made independent of the dynamical parameter by a gauge transformation. 
Computing this constant twist and the corresponding R-matrix in 
the vector representation in the case $\g=sl(n)$, we obtain an explicit
solution of the quantum Yang-Baxter equation on an n-dimensional space. 
This solution is exactly the same as suggested earlier
in \cite{Sch2} on the basis of computer calculations.
The second author conjectured in \cite{Sch2} that it coincides
with the GGS solution, and checked it for $n\le 12$ using a
computer. This conjecture is also proved in \cite{Sch2}
for many special types of
Belavin-Drinfeld triples. 

\begin{rem} \label{r1}
The new dynamical 
twists constructed in this paper
give rise to new examples of Hopf algebroids (quantum groupoids). 
There are two methods to construct a Hopf algebroid 
out of a dynamical twist -- the method of \cite{Xu} and the method
of \cite{EV}, and the results they give are, essentially,
dual to each other.   
For nilpotent triples, when the twist 
is constant (in which case both methods are classical),
the first method yields a new quantum deformation of $U(\g)$
as a Hopf algebra, and the second method yields a 
quantum deformation of the function algebra 
$\text{Fun}(G)$ on the corresponding group. 
We expect that the study of the structure and 
(co)representation theory of these new quantum groups and 
groupoids is a very interesting and fruitful direction of future research.
\end{rem}

\begin{rem} \label{r2}
The results of this paper can be generalized to Kac-Moody
algebras, which we plan to do in a forthcoming paper. 
This generalization is especially interesting in the 
case of affine Lie algebras. In this case, projecting 
the obtained quantum dynamical R-matrices to finite dimensional 
representations, one obtains the quantization of classical dynamical
r-matrices with spectral parameters which are discussed
 in
\cite{ES}. In the case of the Felder and the Belavin r-matrix, it 
was done in \cite{JKOS}, by a method essentially the same as ours. 
\end{rem}

\begin{rem} \label{r3} In \cite{ER}, the authors defined 
the notion of a triangular twist, in order to represent 
quantum determinants of quantum groups as products of 
quasideterminants. 
We would like to emphasize that the twist we constructed 
in this paper is triangular according to the definition 1.2 in \cite{ER}. 
Therefore, the main theorem of \cite{ER}, which claims that the
quantum determinant is a product of quasideterminants, is valid for 
any Belavin-Drinfeld triple, as was anticipated in \cite{ER}. 
\end{rem}

\subsection{Contents}

In Section \ref{bds} we discuss the theory of generalized
Belavin-Drinfeld triples and the corresponding dynamical
r-matrices. 

In Section \ref{dts} we give the construction of the twist ${\mathcal J_T}$
corresponding to a Belavin-Drinfeld triple, 
state the main result (that it satisfies the dynamical 2-cocycle
condition), and show,
by computing the quasiclassical limit of ${\mathcal J_T}$,  
that it defines a quantization of the triple.

In Section \ref{mts} we prove the 2-cocycle condition for ${\mathcal J_T}$.

In Section \ref{gts}, we quantize gauge transformations for dynamical
r-matrices, thus giving an explicit quantization
for all the r-matrices considered. 

In Section \ref{ns}, we discuss the nilpotent case, in which 
the twist can be made non-dynamical by a gauge transformation. 
We compute the explicit form of this twist, which confirms 
a conjecture of the second author \cite{Sch2}.

In Section \ref{es}, we consider examples of quantization of
Belavin-Drinfeld triples.

\subsection{Acknowledgements}
 
The work of P.E. was
partially supported by the NSF grant 9700477, and was partly done
when P.E. was an employee
of the Clay Mathematical Institute as a CMI prize fellow.
O.S. performed this research in part for the Clay Mathematical Institute.
The work of T.S. was supported by the Harvard College Research Program. 
O.S. is grateful to the
Harvard and MIT mathematics departments for hospitality. 
The authors thank Gerstenhaber, Giaquinto, and Hodges for many
useful discussions. 

\section{Generalized Belavin-Drinfeld triples and solutions 
of the classical dynamical Yang-Baxter equation} \label{bds}

\subsection{Notation on semisimple Lie algebras}  Let $\g$ be a
semisimple Lie algebra over $\C$, $\h \subset \g$ a Cartan subalgebra.
Let $\Gamma\subset \h^*$ be a set of simple roots of $\g$.  Fix a
nondegenerate invariant inner product on $\g$ which is positive
definite on the real span of roots.  We can identify $\h$ with $\h^*$
using the inner product.  Let us denote by $h_\alpha$ the image of a
root $\alpha$ under this identification.

Let $\n_\pm$ be the positive and negative nilpotent subalgebras 
of $\g$. We have a decomposition
$\n_\pm=\oplus_{\pm\alpha>0}\g_\alpha$, where $\g_\alpha$ are 
root subspaces of $\g$. Let $e_\alpha$ be a generator 
of the root subspace $\g_\alpha$ for any $\alpha>0$, and let
$f_\alpha$
be a generator of $\g_{-\alpha}$ such that $(e_\alpha,f_\alpha)=1$. 

\subsection{Generalized Belavin-Drinfeld triples and 
dynamical r-matrices}
\begin{defe} \cite{BD}, \cite{S} A generalized Belavin-Drinfeld triple 
$(\Gamma_1, \Gamma_2, T)$ for $\g$ consists of subsets
$\Gamma_1, \Gamma_2 \subset \Gamma$ together with a bijection $T:
\Gamma_1 \rightarrow \Gamma_2$ which preserves the inner product.  
If, in addition,
$\forall \alpha \in \Gamma_1, \exists k \in
\N$ such that $T^k \alpha \notin \Gamma_1$, the triple is said to
be a {\it nilpotent} or {\it ordinary Belavin-Drinfeld triple}.
\end{defe}
  
Given a generalized Belavin-Drinfeld 
triple, 
we define a Lie algebra homomorphism \linebreak $T:\n_+\to \n_+$, by setting
on simple root elements:
$T(e_{\alpha})=e_{T\alpha}$ if $\alpha\in \Gamma_1$, 
and zero otherwise. It is easy to see that such a homomorphism is
well defined. 

For any generalized Belavin-Drinfeld triple, 
let $\l \subset \h$ be the subspace defined by
$\l = \sp(\alpha - T \alpha \mid
\alpha \in \Gamma_1)^\perp$.

It is clear that 
$\l$ is a nondegenerate subspace for the inner product, since
the inner product is positive definite on the real span of roots.

Let us define a useful linear operator 
on the orthogonal complement $\l^\perp$
to $\l$ in $\h$. 
To do this, observe that  
for any $x\in \l^\perp$, there exists a unique $y\in \l^\perp$ such that 
for all 
$\alpha\in \Gamma_1$ one has $(\alpha-T\alpha,y)=(\alpha+T\alpha,x)$. 
It is clear that $y$ depends linearly on $x$. We will write 
$y=C_Tx$. It is easy to check that the operator 
$C_T:\l^\perp\to \l^\perp$ is skewsymmetric.  
This operator is called  the Cayley transform of $T$.

To any   generalized Belavin-Drinfeld
triple $(\Gamma_1, \Gamma_2, T)$, one may 
associate a solution of
the classical dynamical Yang-Baxter equation as follows.  

For a vector space $V$ with a nondegenerate inner product, denote
by $\Omega_V$ the element of $S^2V$ which is inverse to the 
inner
product of $V$. 

Define the function $r_T: \l^* \to (\g \otimes \g)^\l$
by
\begin{equation}
r_T(\lambda)= r-
\frac{1}{2}(C_T\otimes 1)\Omega_{\l^\perp} + a(\lambda)-a^{21}(\lambda)
\end{equation}
where $r:=\frac{1}{2}\Omega_\h +\sum_{\alpha}
e_\alpha  \otimes f_\alpha $ is the standard 
Drinfeld r-matrix, and 
\begin{equation}
a(\lambda)=\sum_{\alpha} \sum_{l=1}^\infty
e^{-l(\lambda,\alpha)} T^l(e_\alpha ) \otimes f_\alpha 
\end{equation}

\begin{prop}[\cite{S}]\label{P:S1} The function $r_T(\lambda)$
is a solution of the classical dynamical Yang-Baxter equation
\begin{equation}
\begin{split}
\sum_i &\left(x_i^{(1)}\frac{\partial}{\partial
x_i}r^{23}(\lambda)-x_i^{(2)}
\frac{\partial}{\partial x_i}r^{13}(\lambda)+x_i^{(3)}\frac{\partial}{\partial
x_i}r^{12}(\lambda)\right)\\
&+[r^{12}(\lambda),r^{13}(\lambda)]+[r^{13}(\lambda),r^{23}(\lambda)]+
[r^{12}(\lambda),r^{23}(\lambda)]=0.
\end{split}
\end{equation}
\end{prop}

\begin{rem} In the expression for
$r_T(\lambda)$, the sum
$\sum_{l=1}^\infty 
e^{-l(\lambda,\alpha)}T^l(e_\alpha ) \otimes f_\alpha $ is finite
if $T$ acts nilpotently on $\alpha$,  and is 
an infinite series convergent to 
a rational function of $e^{(\lambda,\alpha)}$ if
a suitable power of $T$ preserves $\alpha$. 
\end{rem}

This proposition has a straightforward generalization 
to the case of semisimple Lie algebras. 

\subsection{The classification of dynamical r-matrices} 
It is clear that $r_T+r_T^{21}=\Omega_\g$. 
Conversely, it was shown in \cite{S} that if 
$\g$ is simple
then solutions $r_T(\lambda)$ exhaust all solutions of the
classical dynamical Yang-Baxter equation
with this property, up to isomorphism and gauge
transformations. More precisely, the result is as follows. 

Let $H$ be the Cartan subgroup of the Lie group $G$ corresponding to
 $\g$, whose Lie algebra is $\h$.  By a gauge transformation of 
$\mathbf r: \l^*\to (\g \o \g)^\l$,
 we mean a transformation 
\begin{equation}
\mathbf r \mapsto \mathbf r^g:=
(\Ad g \o \Ad g) 
(\mathbf r - (g^{-1}dg) + (g^{-1}dg)^{21})
\end{equation} 
for $g: \l^*
\rightarrow H$ a meromorphic function
(here the differential 1-form $g^{-1}dg$ on $\l^*$ with values in $\h$ is
 regarded as a function on $\l^*$ with values in $\l\otimes \h$).
Now, we have the following theorem:

\begin{thm} \cite{S} Let $\g$ be simple, and let 
$\l' \subset \h$ be any subalgebra which is nondegenerate with 
respect to the inner product on $\h$.
Let $\mathbf r:{\l'}^*\to (\g \o \g)^{\l'}$
be any solution of the classical dynamical 
Yang-Baxter equation satisfying $\mathbf r + \mathbf r^{21}=\Omega_\g$.
Then, $\mathbf r$ may be obtained from $r_T(\lambda)$ 
for a suitable triple by an
automorphism of 
$\g$ and a gauge transformation.
\end{thm}

The Belavin-Drinfeld result for the ordinary classical
Yang-Baxter
equation (CYBE)
\cite{BD} is easily obtained from the
above classification.  Namely, any solution $\mathbf r \in \g \o \g$ to the
CYBE such that $\mathbf r+\mathbf r^{21}=\Omega_\g$ 
is equivalent under an automorphism of $\g$ to a solution of the
form 
\begin{equation}
r_{T,s}= r
-s 
+\sum_{\alpha} \sum_{l=1}^\infty
T^l(e_\alpha ) \wedge f_\alpha 
\end{equation}
 for an ordinary
(nilpotent) Belavin-Drinfeld triple $(\Gamma_1, \Gamma_2, T)$,
where $s\in \Lambda^2\h$ is a solution of the equations 
\begin{equation} \label{r0h}
((\alpha - T
\alpha) \o 1) s = \frac{1}{2}((\alpha + T \alpha) \o 1)
\Omega_\h, \forall \alpha \in \Gamma_1.
\end{equation}
This solution $r_{T,s}$ 
(which is independent of $\lambda\in \l^*$) 
can be obtained by a gauge  transformation from 
$r_T(\lambda)$. 

\begin{rem}
It is not difficult to generalize the results of this section to the semisimple
case. 
\end{rem}

\section{The dynamical twist $\J_T$} \label{dts}

\subsection{The extension of $T$ to an orthogonal automorphism of
${\mathfrak l}$} \label{te}

The goal of Section \ref{dts} is to construct a twist which provides a
quantization of the classical dynamical r-matrices of Section \ref{bds}. In order
to define this twist, we need to extend the map $T$ to a linear map from 
$\h$ to $\h$.

Fix a finite dimensional semisimple Lie algebra 
$\g$ and a generalized Belavin-Drinfeld triple $(\Gamma_1, \Gamma_2,
T)$.  Let $\Gamma_3 \subset \Gamma_1$ be the largest
$T$-invariant subset of $\Gamma_1$.  
Set $\h_i = \text{Span}(\Gamma_i)$ for $i \in
\{1,2,3\}$. It is clear that $T$ extends 
naturally to a linear map $T: \h_1 \rightarrow
\h_2$, which we will also denote by $T$. 

\begin{lemma}  One has
 $\h_1 \cap \l = \h_3^{ T}$. Moreover, $\h_1 + \l = \h$.
\end{lemma}

{\it Proof.}  
Let $x \in \h_1 \cap \l$ be a real element.    
Then $(T x, T x) = (x,x) =  (x, T x) = (T x, x)$. Thus, $(x-T x, x -
T x) = 0$, 
so the positive definiteness of the form on real elements implies $x = T x$. 
Since $ T$ is nilpotent as a
map 
$\Gamma_1 \setminus \Gamma_3 \rightarrow \Gamma_2 \setminus
\Gamma_3$, we get $x\in \h_3^T$. This proves the first
statement. 

To prove the second statement, it is enough to notice that the
 rank of the system of linear equations 
$\alpha(y)=T\alpha(y)$, $\alpha\in\Gamma_1$ (with respect to $y\in\h$), is
 $|\Gamma_1|-|\Gamma_3/T|$. Thus, $dim({\mathfrak
 l})=|\Gamma|-|\Gamma_1|
+|\Gamma_3/T|$. So the second statement follows from the first 
statement. \qs

\begin{corr}  There is a unique extension $T: \h \rightarrow \h$ of 
$ T:\h_1\to\h_2$ which is equal to the identity on
$\l$. This extension is an orthogonal operator.
\end{corr}

{\it Proof.}  By the lemma, $\h = \h_1 + \l$ and $T$ is the identity
on $\h_1 \cap \l$, which immediately implies that $T$ admits a
unique extension as desired. The fact that the obtained extension
  is orthogonal is straightforward to verify. \qs

\subsection{The construction of the twist}
 
Let $q=e^{\hbar/2}$, where $\hbar$ is a formal parameter, and 
let $U_q(\g)$ be the Drinfeld-Jimbo quantum universal enveloping 
algebra, which is a 
quantization of the Lie bialgebra 
$(\g,r)$ (see \cite{CP}, p 281). Let $\R=1+\hbar r+...$ 
be its universal R-matrix. We choose the coproduct of $U_q(\g)$
in such a way that the ``Cartan part'' of the R-matrix is
$q^{\Omega_\h}$. 

Our goal in this subsection
is to introduce a dynamical twist $\J_T(\lambda)$, $\lambda\in \l^*$ 
in the (completed) tensor square of $U_q(\g)$ 
(i.e. a solution of the dynamical 2-cocycle condition) which would twist the
universal R-matrix of $U_q(\g)$ to a matrix $R_T = 1 +
\hbar r_T
 + O(\hbar^2)$.  

Define the degree of an element of $U_q(\g)$ by $deg(e_i)=1$,
$deg(f_i)=-1$, $deg(h_i)=0$, and $deg(xy) = deg(x)+deg(y)$.

Let $U_+$
denote the Hopf subalgebra of $U_q(\g)$ generated by 
elements $h_\alpha,e_i$ and let 
$U_-$
denote the Hopf subalgebra of $U_q(\g)$ generated by 
elements $h_\alpha,f_i$.

Define Hopf algebra homomorphisms
$T:U_+\to U_+$, $T^{-1}:U_-\to U_-$ as follows: 
$T^{\pm 1}$ on $\h$ is defined as in Section \ref{te},  
$T(e_\alpha)=e_{T\alpha}$ if $\alpha\in \Gamma_1$, 
$T(e_\alpha)=0$ for other simple roots $\alpha$, 
$T^{-1}(f_\alpha)=f_{T^{-1}\alpha}$ if 
$\alpha\in \Gamma_2$, 
$T^{-1}(f_\alpha)=0$ on other simple roots. It is easy to show 
that such homomorphisms exist and are unique. 

Set $Z = (\frac{1}{2}(C_T-1) \o 1) \Omega_{\l^\perp}$.

Let $W_2$ be the vector space of rational functions of 
$e^{(\lambda,\alpha)}$ (for simple roots $\alpha$)
with values in $(U_+\otimes U_-)^{\l}$, 
where $\otimes$ is the completed tensor product
in the $\hbar$-adic topology. 

\begin{rem}
We use the notation $W_2$ since this space consists of
2-component tensors. 
\end{rem}

Define a linear operator $A_L^2:W_2\to W_2$ as follows:

\begin{gather} \label{a2l}
A_L^2 X = (T \Ad e^{-\lambda} \o 1)(\R Xq^{-\Omega_\l})
\end{gather}

\begin{rem} The notation $A_L^2$ has the following motivation:
this is an operator on 2-component tensors, which applies $T$ to
the left component (so $L$ stands for ``left''). 
\end{rem}

Let $I_\pm$ denote the kernels of the projections of $U_\pm$
to elements of zero degree. 

The construction of the dynamical twist depends on the following
proposition. 

\begin{prop} \label{jp} There exists a unique element $\J_T \in W_2$ such that

1) $\J_T - q^Z \in I_+ \o
I_-$,

2) $\J_T$ satisfies the ``modified ABRR equation,'' $A_L^2 \J_T = \J_T$.
\end{prop}

\begin{rem} If $T=id$, this proposition is contained in the
paper 
\cite{ABRR}, which motivates the terminology ``the modified ABRR equation''.
\end{rem}

{\it Proof.}  
The statement is equivalent to the claim that there exists a
unique $X_0$ belonging to $1+(I_+\otimes I_-)^\l$ such that 
\begin{equation}
X_0=(T\Ad e^{-\lambda}\otimes 1)(\R X_0 q^{-\Omega_\h})
\end{equation}
(then $\J_T=X_0 q^Z$). Let us write $X_0$ as $1+\sum_{j\ge 1}
X_0^j$,
where $X_0^j$ are the
terms of degree $j$ in the first component. Then 
the above equation can be written as a system of equations
labeled by degree $j\ge 1$:
\begin{equation}
X_0^j=(T\Ad e^{-\lambda}\otimes
1)(q^{\Omega_\h}X_0^jq^{-\Omega_\h})+...
\end{equation}
where $...$ stands for terms that involve $X_0^i$ for $i<j$. 
(The zero degree equation is obviously satisfied, so we don't
need to include it).
It is obvious that the operator 
\begin{equation}
1-(T\Ad e^{-\lambda}\otimes 1)\Ad q^{\Omega_\h}
\end{equation} 
on $W_2$ is invertible 
for generic $\lambda$. Thus, the last equation admits a unique 
solution for all $j$, which allows one to compute $X_0^j$
recursively. The proposition is proved. 
\qs

\subsection{The main theorem}

The main theorem of this paper is the following:

\begin{thm} \label{mt}
The element $\J_T$ satisfies the dynamical cocycle condition,
\begin{equation} \label{dcoc}
\J_T^{12,3}(\lambda) \J_T^{12}(\lambda-\frac{1}{2}\hbar h^{(3)}) =
\J_T^{1,23}(\lambda) \J_T^{23}(\lambda + \frac{1}{2} \hbar h^{(1)}).
\end{equation}
\end{thm}

Here, by $\J_T^{12,3}$ we mean $(\Delta \o 1) (\J_T)$ where $\Delta:
U_q(\g) \rightarrow U_q(\g) \o U_q(\g)$ is the coproduct of $U_q(\g)$
and similarly $\J_T^{1,23} = (1 \o \Delta)(\J_T)$, and
$\lambda-\frac{1}{2}\hbar h^{(3)}$ is defined as follows. If
$y_1,..,y_r$ is a basis of $\l$, and
$\lambda=(\lambda^1,..,\lambda^r)$ is the coordinate representation of
$\lambda$ with respect to this basis, then $\lambda-\frac{1}{2}\hbar
h^{(3)}= (\lambda^1-\frac{1}{2}\hbar
y_1^{(3)},...,\lambda^r-\frac{1}{2}\hbar y_r^{(3)})$, and for
any meromorphic function $S(\lambda)$ we set
$S(\lambda-\frac{1}{2}\hbar h^{(3)}) = S(\lambda) -\frac{1}{2} \hbar
\sum_i \frac{\partial S}{\partial y_i}\bigl(\lambda\bigr) h_i +
\ldots$ (the Taylor expansion). The definition
of $\lambda + \frac{1}{2} \hbar h^{(1)}$ is similar.  

The proof of the theorem is given in Section \ref{mts}.

Now define 
\begin{equation}
R_T(\lambda):=(\J_T^{21})^{-1}(\lambda) \R \J_T(\lambda).
\end{equation} 

\begin{corr} The function $R_T(\lambda)$ satisfies the
  (symmetrized) quantum dynamical Yang-Baxter equation
\begin{multline}
R_T^{12}(\lambda+\frac{1}{2}\hbar h^{(3)})
R_T^{13}(\lambda-\frac{1}{2}\hbar h^{(2)})
R_T^{23}(\lambda+\frac{1}{2}\hbar h^{(1)}) \\ =
R_T^{23}(\lambda-\frac{1}{2}\hbar h^{(1)})
R_T^{13}(\lambda+\frac{1}{2}\hbar h^{(2)})
R_T^{12}(\lambda-\frac{1}{2}\hbar h^{(3)}).
\end{multline}
\end{corr} 

The proof of the corollary is straightforward using the main
theorem. 

\begin{rem} The symmetrized quantum dynamical Yang-Baxter
equation was first considered by Felder. 
It is equivalent to the
nonsymmetrized quantum dynamical Yang-Baxter 
equation 
\begin{equation}
R^{12}(\lambda-\hbar h^{(3)})
R^{13}(\lambda)
R^{23}(\lambda-\hbar h^{(1)})=
R^{23}(\lambda)
R^{13}(\lambda-\hbar h^{(2)})
R^{12}(\lambda).
\end{equation}
considered by many authors, by the change of variable 
$\lambda=-(\lambda'+\frac{1}{2}(h^{(1)}+h^{(2)}))$. 
\end{rem}

\subsection{The quasiclassical limit}
\begin{prop} One has $R_T(\lambda)= 1 + \hbar r_T(\lambda) \pmod{\hbar^2}$,
where $r_T(\lambda)$ is as given in the previous section.
In other words, $R_T(\lambda)$ is a quantization of
$r_T(\lambda)$. 
\end{prop}

{\it Proof.}  The quasiclassical limit 
of the modified ABRR equation has the form
\begin{equation}
x=(T\Ad e^{-\lambda}\otimes 1)(r+x-\frac{1}{2}\Omega_\l),
\end{equation}
where $x\in (\g\otimes \g)^\l$ has Cartan part $Z$. Solving this equation,
we obtain that $\J_T= 1 + \hbar (\frac{Z}{2} + a)
\pmod{\hbar^2}$. Since $\R=1+\hbar r\pmod{\hbar^2}$, the proposition
is proved.  \qs

\subsection{The nilpotent case}
If the Belavin-Drinfeld triple is
nilpotent, the element $\J_T$ can be written as a finite
product, as follows.   

For a nilpotent triple, let $n$ be the
largest integer so that $T^n$ is defined on some $\alpha \in
\Gamma_1$.
 
\begin{prop} \label{nj} In the nilpotent case, the solution 
$\J_T$ of the modified ABRR equation can be written in the form
\begin{equation}
\J_T(\lambda) =(T\Ad e^{-\lambda}\o 1)(\R)...(T^n\Ad
e^{-n\lambda}\o 1)(\R) 
q^{Z-[(T+...+T^n)\otimes 1](\Omega_\h)}.
\end{equation}
\end{prop}

The proof of this proposition 
is obtained by substituting the modified ABRR equation
into itself $n$ times. 

\subsection{The infinite product formula for $\J_T$}

A product formula for $\J_T$ similar to the above exists for an
arbitrary
(not necessarily nilpotent) $T$, but in the non-nilpotent case 
the product is infinite. To write down the general formula, let 
$\R=\R_0q^{\Omega_\h}$, and 
\begin{equation}
\R_0^m(\lambda)=\Ad q^{([T+...+T^{m-1}]\o 1)\Omega_\h}
(T^m\Ad e^{-m\lambda}\otimes
1)\R_0.
\end{equation} 
Then we have 
\begin{equation}
\J_T(\lambda)=[\prod_{m=1}^\infty \R_0^i(\lambda)]q^Z.
\end{equation}
This product is clearly convergent in the topology of formal 
power series in $e^{-(\lambda,\alpha)}$.
If $T$ is nilpotent, the product becomes finite and we get the
formula
from the previous section. 

\section{Proof of the main theorem} \label{mts}

\subsection{The right component version of the modified ABRR
equation}

To prove Theorem \ref{mt} we introduce a ``right-component version''
of the modified ABRR equation. Define the linear operator 
$A_R^2:W_2\to W_2$ by
\begin{equation} \label{ar2}
A_R^2 X = (1\otimes T^{-1} \Ad e^{\lambda})(\R Xq^{-\Omega_\l})
\end{equation}

\begin{lemma} The operators $A_L^2, A_R^2$ commute. \end{lemma}
{\it Proof.} This follows immediately from the fact that $(T \o 1)
(\R) = (1 \o T^{-1}) (\R)$, which is true because 
$T$ preserves the inner product on $\h$. \qs

\begin{corr} $\J_T$ is the unique solution to 
the system of equations $A_R^2 X = X, A_L^2 X = X$ with
$\J_T-q^Z\in (I_+\otimes I_-)^\l$.
\end{corr}

{\it Proof.}  We have $A_L^2 A_R^2 \J_T = A_R^2 A_L^2 \J_T =
A_R^2 \J_T$ 
so that
$A_R^2 \J_T$ and $\J_T$ are both solutions to $A_L^2 X = X$
with zero degree term $q^Z$. Hence $\J_T = A_R^2 \J_T$
by Proposition \ref{jp}. \qs

\subsection{The 3-component versions of the modified ABRR equation}

Now, we introduce ``3-component versions'' of the modified ABRR
equation.
Let $W_3$ be 
the vector space of 
rational functions of $e^{(\lambda,\alpha)}$ for simple roots
$\alpha$, with values in $(U_+\otimes U\otimes U_-)^\l$,
where $\otimes$ is the completed tensor product
in the $\hbar$-adic topology. 

Define linear operators $A_L^3,A_R^3:W_3\to W_3$ by 
\begin{gather}
A_L^3 X = (T \Ad e^{-\lambda} \o 1 \o 1)(\R^{13} \R^{12}
Xq^{-\Omega_\l^{13}-
\Omega_\l^{12}}), \\ A_R^3 X = (1 \o 1 \o T^{-1}
\Ad e^{\lambda})(\R^{13} \R^{23} X
q^{-\Omega_\l^{13}-
\Omega_\l^{23}}) 
\end{gather}

The 3-component versions of the modified ABRR equation are 
$A_L^3X=X,A_R^3X=X$. They are obtained from the
2-component equations by comultiplication of the component in
which there is no action of $T$. 

\begin{lemma} The operators $A_L^3$ and $A_R^3$ commute. \end{lemma}
{\it Proof.}  This reduces to showing that 
\begin{gather}
(T \Ad e^{- \lambda} \o 1 \o 1)
(\R^{13} \R^{12} (1 \o 1 \o T^{-1} \Ad e^{\lambda})(\R^{13}
\R^{23}) )=
\\ (1 \o 1 \o T^{-1} \Ad e^{\lambda})(\R^{13} \R^{23} (T \Ad 
e^{-\lambda} \o 1 \o 1) (\R^{13} \R^{12})).
\end{gather}
 Let $\bar \R = (T
\Ad e^{-\lambda} \o 1) \R = (1 \o T^{-1} \Ad e^{\lambda}) \R$
and $ \tilde \R = (T \Ad e^{-\lambda} \o T^{-1} \Ad e^{\lambda})
\R$.  Then we need to check 
$\bar \R^{13} \bar \R^{12} \tilde \R^{13} \bar \R^{23}
= \bar \R^{13} \bar \R^{23} \tilde \R^{13} \bar \R^{12}$, which
follows (after cancelling the first factor) from the
quantum Yang-Baxter equation
 for $\R$, applying $T \Ad e^{-\lambda}$ in the first component and
$T^{-1} \Ad e^{\lambda}$ in the third one.\qs

\begin{lemma} If there exists a solution $X$ of $A_L^3 X = A_R^3 X = X$
such that \\ $X - q^{Z_{12} + Z_{13} + Z_{23}} \in I_+ \o U(\g) \o U(\g) +
U(\g) \o U(\g) \o I_-$, it is unique.
\end{lemma}

{\it Proof.}  It is enough to show 
that such a solution $X$ is unique for the equation 
$A_L^3A_R^3X=X$. 
Let us make a change of variable $X=X_0q^{Z_{12}+Z_{13}+Z_{23}}$,
and write $X_0=1+\sum_{k,l\ge 0: k+l>0}X_0^{k,l}$, where 
$X_0^{k,l}$ is the part of $X_0$ having degree
$k$ in the first component and $-l$ in the third component.
 It is easy to check that the
equation for $A_L^3A_R^3X=X$ transforms to the system of equations
\begin{equation}
X_0^{k,l}=(T\Ad e^{-\lambda}\otimes 1\otimes T^{-1}\Ad
e^{\lambda})(q^WX_0^{k,l}q^{-W})+..., k+l>0
\end{equation}
where $W=\Omega_\h^{12}+\Omega_\h^{23}+[1\otimes
(1+T)]\Omega_\h^{13}$, 
and $...$   
stands for terms that involve $X_0^{k',l'}$ with $k'+l'<k+l$ 
(The zero degree equation is obviously satisfied, so we don't
need to include it).
It is obvious that the operator 
\begin{equation}
1-(T\Ad e^{-\lambda}\otimes 1\otimes T^{-1} \Ad e^{\lambda})
\Ad q^{W}
\end{equation} on $W_3$ is invertible 
for generic $\lambda$. Thus, the last equation admits a unique 
solution for all $k,l$, which allows one to compute $X_0^{k,l}$
recursively. The lemma is proved. 
\qs

Now we complete the proof of the main theorem. 
It is obvious that $\J_T^{12,3}(\lambda)
\J_T^{12}(\lambda-\frac{1}{2}\hbar h^{(3)})$ is a solution of $A_R^3X=X$ and
that
$\J_T^{1,23}(\lambda) \J_T^{23}(\lambda + \frac{1}{2}\hbar h^{(1)})$ is a
solution of $A_L^3X=X$.  So, by virtue of the previous lemma, 
to prove the main theorem it is sufficient to prove the
following. 

\begin{lemma} (i) $X = \J_T^{12,3}(\lambda) \J_T^{12}(\lambda
  -\frac{1}{2} 
\hbar h^{(3)})$ is a solution of $A_L^3 X = X$.

(ii) $X = \J_T^{1,23}(\lambda) \J_T^{23}(\lambda + \frac{1}{2}
\hbar h^{(1)})$ is
a solution of $A_R^3 X = X$.
\end{lemma}

{\it Proof.} (i) As we have mentioned, the element 
$X = \J_T^{12,3}(\lambda) \J_T^{12}(\lambda -
\frac{1}{2} \hbar h^{(3)})$ satisfies $A_R^3 X = X$.  Since $A_L^3$ and
$A_R^3$ commute, the element $Y = A_L^3 X$ is also a solution of $A_R^3 X = X$.
Since any element of $W_3$ invariant under
$A_R^3$ is uniquely determined by
its part of zero degree in the third component, it suffices to
show that $A_L^3 X$ has the same part of zero degree in the third
component as $X$.  Call the former $Y_0$ and the latter $X_0$.
Clearly $X_0 = q^{Z^{13} + Z^{23}} \J_T^{12}(\lambda - \frac{1}{2}
\hbar h^{(3)})$.  Now, we find from $Y=A_L^3X$ that
\begin{equation}
Y_0 = (T\Ad e^{-\lambda} \o 1 \o 1)(q^{\Omega_\h^{13}}
\R^{12} q^{Z^{13} + Z^{23}} \J_T^{12}(\lambda - \frac{1}{2} \hbar
h^{(3)}))
q^{-\Omega_\l^{12} - \Omega_\l^{13}}.
\end{equation}
  Since $[Z^{13} + Z^{23}, \R^{12}]
= 0$, we may rewrite this as 
\begin{equation}
Y_0 = q^{Z^{23}} (T\Ad e^{-\lambda} \o
1 \o 1)(q^{\Omega_{\l^\perp}^{13}+\Omega_{\l}^{13}+Z^{13}} \R^{12}
\J_T^{12}(\lambda -\frac{1}{2} \hbar h^{(3)}))
q^{-\Omega_\l^{12}-\Omega_\l^{13}}.
\end{equation}
  Since $(T \o 1)(Z +
\Omega_{\l^\perp}) = Z$, we have 
\begin{equation}
Y_0=q^{Z^{23}+ Z^{13}} (T\Ad e^{-\lambda} \o
1 \o 1)(\Ad q^{\Omega_\l^{13}} (\R^{12} \J_T^{12}(\lambda - \frac{1}{2}\hbar
h^{(3)}))) q^{-\Omega_\l^{12}}.
\end{equation}
  Next, note that $\Ad q^{h^{(3)}} \o 1
\o 1 = \Ad q^{\Omega_\l^{13}}$, so that 
\begin{equation}
Y_0 = q^{Z^{23}+ Z^{13}}(T \Ad
e^{- \lambda + \frac{1}{2}\hbar 
h^{(3)}} \o 1 \o 1) (\R^{12} \J_T^{12}(\lambda -
\frac{1}{2} \hbar h^{(3)})) q^{-\Omega_\l^{12}}.
\end{equation}
Changing $\lambda$ to $\lambda-\frac{1}{2}\hbar h^{(3)}$ in the
modified ABRR equation for $\J_T$, we see from the last equation
that 
$Y_0=X_0$, as desired.   

(ii)  This is proved analogously to (i). \qs

The main theorem is proved. 

\section{Quantization of gauge transformations} \label{gts}

In this section we provide a quantization for all gauge
transformations, which yields a quantization for all dynamical
r-matrices considered in Section \ref{bds}. 

Let $g:\l^*\to H$ be a meromorphic function. 

\begin{lemma} If $\J(\lambda)\in U_q(\g)\otimes U_q(\g)$ 
is a solution of the dynamical cocycle condition of Theorem \ref{mt}, 
then so is 
\begin{equation}
\J^g(\lambda)=(g(\lambda)\otimes g(\lambda))
\J(\lambda)(g^{-1}(\lambda-\frac{1}{2}\hbar h^{(2)})\otimes 
g^{-1}(\lambda+\frac{1}{2}\hbar h^{(1)})).
\end{equation}
\end{lemma}

The proof of this lemma is straightforward. 

\begin{corr}
The element 
\begin{equation}
\J_T^g(\lambda)=(g(\lambda)\otimes g(\lambda))
\J_T(\lambda)(g^{-1}(\lambda-\frac{1}{2}\hbar h^{(2)})\otimes 
g^{-1}(\lambda+\frac{1}{2}\hbar h^{(1)}))
\end{equation} 
satisfies the dynamical 2-cocycle condition. 
The element $R_T^g(\lambda)=(\J_T^{g})^{21}(\lambda)^{-1}\R
\J^g_T(\lambda)$ satisfies the quantum dynamical Yang-Baxter 
equation, and is a quantization of the solution $r_T^g(\lambda)$
of the classical dynamical Yang-Baxter equation, which is
obtained from $r_T$ by the gauge transformation $g$. 
 
\end{corr}

The proof of the corollary is by an easy direct calculation. 

\section{The nilpotent case} \label{ns}

\subsection{The $\lambda$-independent twist}

In the case of nilpotent (or ordinary)
Belavin-Drinfeld triples, the dynamical twist $\J_T(\lambda)$ may be
transformed by a gauge transformation
into an ordinary twist that does not depend on $\lambda$,
and hence satisfies the ordinary
(non-dynamical) 2-cocycle condition 
\begin{equation}
\J^{12,3} \J^{12} = \J^{1,23} \J^{23}.
\end{equation}  
This yields an explicit quantization for all 
non-dynamical r-matrices defined in Section \ref{bds}.

Namely, in the setting of Section \ref{gts}, set
$g(\lambda)=e^{Q\lambda}$, where $Q;\l^*\to \h$ is a linear map. 
We have $\J_T^g(\lambda)=(\Ad e^{Q\lambda}\otimes 
\Ad e^{Q\lambda})\J^T(\lambda)q^{Q^{21}-Q}$, where 
in the last factor we understand $Q$ as an element of $\l\otimes \h$. 
According to Section \ref{gts}, this element satisfies the dynamical
2-cocycle condition. 

Now choose a solution 
$s$ of equation \eqref{r0h} and take $Q$ to be the component of $s$ 
in $\l \o \l^\perp$.  Denote $\J_T^g$ by $\J_{T,s}$. 

\begin{thm}  $\J_{T,s}$ is independent of 
$\lambda$ (i.e. $\J_{T,s}: \l^* \rightarrow
U_q(\g) \o U_q(\g)$ is constant.) 
Hence, $\J_{T,s}$ satisfies the non-dynamical 
2-cocycle condition $\J^{12,3}\J^{12}=\J^{1,23}\J^{23}$.
\end{thm}

{\it Proof.} We may write as in Proposition \ref{nj} that 

\begin{multline} \label{jpt}
\J_{T,s} = (\Ad e^{Q\lambda}\o \Ad e^{Q\lambda}) 
[(T \Ad e^{-\lambda}  \o 1) (\R) \cdots 
(T^n \Ad e^{-n \lambda} \o 1) (\R) \\ q^{Z - Q + Q^{21} -
[(T+...+T^n) \o 1]( \Omega_{\h})}].
\end{multline}

Now, we note that 
\begin{equation}
(\Ad e^{Q\lambda}\o 
\Ad e^{Q\lambda}) (e_{T^l \alpha}\o f_\alpha)
= e^{(Q\lambda,T^l \alpha - \alpha)} e_{T^l \alpha}\o 
f_\alpha.
\end{equation} 
But by equation \eqref{r0h}, 
\begin{equation}
(Q\lambda,T^l \alpha - \alpha)=
(s,\lambda\otimes (T^l\alpha-\alpha))=
\sum_{j=0}^{l-1}\frac{1}{2}(\lambda,(T^j+T^{j+1})(\alpha))=
l(\lambda,\alpha)
\end{equation}
(as $(\lambda,T\alpha)=(\lambda,\alpha)$). 

But it is easy to show 
from the ABRR equation that in the nilpotent case
$\J_T$ is a linear combination 
with constant coefficients of products of $
e^{-l(\lambda,\alpha)}e_{T^l \alpha}\o
f_\alpha$ 
and elements from $U_q(\h)^{\o 2}$. 
Therefore, 
the above calculation precisely implies that there is no
$\lambda$-dependence
 in 
$\J_{T,s}$.  \qs

\begin{corr}
\begin{equation}
\J_{T,s}=(T\otimes 1)(\R)...(T^n\otimes 1)(\R)
q^{-s-\frac{1}{2}\Omega_{\l^\perp}-[(T+...+T^n)\otimes 1](\Omega_\h)}.
\end{equation}
\end{corr}

{\it Proof.} Since there is no $\lambda$-dependence, we 
can set $\lambda$ to $0$ in the expression for $\J_{T,s}$, which 
after a short calculation yields the result. \qs

\begin{corr} The element $R_{T,s}=(\J_{T,s}^{21})^{-1}
\R\J_{T,s}$ is a solution of the quantum Yang-Baxter equation 
$R^{12}R^{13}R^{23}=R^{23}R^{13}R^{12}$ 
which is a quantization of $r_{T,s}$. 
\end{corr} 

The proof is straightforward. 

\begin{rem} The fact that the element $J_{T,s}$ is a twist 
(i.e. satisfies the 2-cocycle condition) can be proved 
without ever mentioning ``dynamical'' objects, 
along the lines of Section \ref{mts}. Namely, 
the new equations would be as before
but with operators $(A_L^2)',(A_R^2)',(A_L^3)',(A_R^3)'$ defined as follows:
\begin{gather}
(A_L^2)' X = (T \o 1) (\R X q^{Q-Q^{21}}) q^{Q^{21} -Q- \Omega_\l}, \\
(A_R^2)' X = (1 \o T^{-1}) (\R X q^{Q-Q^{21}}) q^{Q^{21} - Q-
  \Omega_\l}, 
\\
(A_L^3)' X = (T \o 1 \o 1) (\R^{13} \R^{12} X q^{-Q^{21} - Q^{31} + Q^{12} 
+ Q^{13}}) q^{-Q^{12} - Q^{13} + Q^{21} + Q^{31} - \Omega_\l^{12} - 
\Omega_\l^{13}}, \\
(A_R^3)' X = (1 \o 1 \o T^{-1}) (\R^{13} \R^{23} X q^{-Q^{31} -
  Q^{32} + Q^{13} + Q^{23}}) q^{-Q^{13} - Q^{23} + Q^{31} + Q^{32}
  - \Omega_\l^{13}
 - \Omega_\l^{23}}. 
\end{gather}
These equations are obtained from the corresponding equations of
Section \ref{mts} by the gauge trasformation $e^{Q\lambda}$ as in Section
5, and then setting $\lambda$ to 0. 

One may show that $X = \J_{T,s}$ is the unique solution to
$(A_L^2)' X = X$
with a suitable part of zero degree in the first component, 
and also satisfies $(A_R^2)' X = X$, and that $(\J_{T,s})^{12,3}
(\J_{T,s})^{12}$ and $(\J_{T,s})^{1,23} (\J_{T,s})^{23}$ are both
equal to the
unique 
solution to the
system $(A_L^3)' Y = (A_R^3)' Y = Y$ with a suitable zero degree part.
This implies that $\J_{T,s}$ satisfies the 2-cocycle condition. 
\end{rem}

\subsection{Explicit calculation of $\J_{T,s}$ in the vector 
representation for $\g={\mathfrak {sl}}_n$}

Consider the case $\g = \mathfrak{sl}(n)$. Let $\Gamma =
\{\alpha_1,\ldots,\alpha_{n-1}\}$ be the set of simple roots where
$\alpha_i=v_i-v_{i+1}$, and $v_i$ is the standard basis of 
$\mathbb Z^n$. For a root $\alpha=v_i-v_j$, let $e_\alpha=e_{ij}$ 
be the corresponding elementary matrix. 
 
We have the $n$-dimensional representation, $\phi: U_q(\g) \rightarrow
Mat_n(\C)$, given by $\phi(e_{\alpha_i}) = e_{i,i+1},
\phi(f_{\alpha_i}) = e_{i+1,i}, \phi(h_{\alpha_i}) = e_{ii} -
e_{i+1,i+1}$. Let us calculate explicitly 
the matrix \linebreak $(\phi\otimes \phi)(\J_{T,s})$.

Let $\Gamma_1^{(k)} \subset \Gamma_1$ be the subset of all simple
roots $\alpha$ on which $T^k$ is defined.  Let $\tilde \Gamma_1^{(k)}
\subset \text{Span}(\Gamma_1^{(k)})$ denote the subset of positive
roots.  Whenever $\alpha \in \tilde \Gamma_1^{(k)}$ is not simple, let
$C_{\alpha,k} = 1$ if $T^k$ reverses the orientation of $\alpha$ (as a
segment on the Dynkin diagram) and $0$ otherwise. For simple roots,
let $C_{\alpha,k}$ be 0.

Denote by $|\alpha|$ the number of simple roots in a positive root
$\alpha$. 
For two positive roots $\alpha,\beta$ we will write 
$\alpha\lessdot\beta$ if $\alpha=v_i-v_j$, and 
$\beta=v_j-v_k$. Write $\alpha\prec\beta$ if 
$T^k\alpha=\beta$ for some $k>0$.   

For $\alpha\prec\beta$, define 
\begin{equation}
L_{\alpha,\beta} =
\frac{1}{2} [\alpha \lessdot \beta] - \frac{1}{2}
[\beta \lessdot \alpha]+ 
[\exists \gamma, \alpha \prec \gamma \prec \beta,
\alpha \lessdot \gamma] - [\exists \gamma, \alpha \prec \gamma \prec \beta,
 \gamma \lessdot\alpha],
\end{equation}
where 
$[\text{statement}] = 0$ if
statement is false and otherwise $[\text{statement}] = 1$.

\begin{prop} The action of $\J_{T,s}(\lambda)$ on
the tensor product of two vector representations is given by 
\begin{equation}
(\phi\otimes \phi)(\J_{T,s})=
q^{-\frac{1}{2}\Omega_\h}
J_1\cdots J_n q^{-s+\frac{1}{2}\Omega_{\l}},
\end{equation}
where 
\begin{equation}
J_k=1+\sum_{\alpha \in \tilde
\Gamma_1^{(k)}} (-q)^{(|\alpha| - 1)C_{\alpha,k}}
q^{L_{\alpha,T^k\alpha}} (q-q^{-1}) e_{T^k\alpha} \o
e_{-\alpha}
\end{equation}
\end{prop}

{\it Proof.} The proposition follows by a direct calculation 
from the explicit expression of the R-matrix in \cite{KhT} 
and the above formula for $\J_{T,s}$.\qs
 
Let 
\begin{equation}
R
=q^{-1/n}(q \sum_i e_{ii}\otimes e_{ii}+
\sum_{i\ne j}e_{ii}\otimes e_{jj}+
(q-q^{-1})\sum_{i<j}e_{ij}\otimes e_{ji}).
\end{equation} 
be the standard R-matrix in the vector 
representation (it is a quantization of $r$ in the vector representation). 

\begin{corr} The element 
\begin{equation}
R(T,s)=q^{-s}(J^{21}_n)^{-1}...(J^{21}_1)^{-1}
RJ_1...J_nq^{-s}.
\end{equation} 
of $\End(\C^n)\otimes \End(\C^n)$
satisfies the quantum Yang-Baxter equation, and is a quantization of 
$r_T$ in the vector representation. 
\end{corr}

The proof is clear, noting that, in the vector representation, 
$\Omega_\l$ commutes with anything that is invariant under $\l$. 

This corollary proves Part 1 of
conjecture 1.2 in \cite{Sch2} since the element
$R_J$ considered there is just $q^{\frac{1}{n}}R(T,s)^{21}$.  The
element $R_J$ coincides with the GGS R-matrix in all checked cases,
as detailed in \cite{Sch2}.

\section{Examples of the twist $\J_T$} \label{es}

In this section we compute the twist $\J_T$ for $\g = \mathfrak{sl}(n)$,
evaluated in the representation $\phi: U_q(\g) \rightarrow Mat_n(\C)$
in two particularly simple cases, when 
$\Gamma=\Gamma_1=\Gamma_2$.  In particular, we give the twist for
all $n$ where $T = id$
(this is the case considered in \cite{ABRR},\cite{JKOS}), 
and for $\g = \mathfrak{sl}(3)$ where
$T(\alpha_1) = \alpha_2, T(\alpha_2) = \alpha_1$ (the ``flip'' map).

\begin{prop}  For the triple $(\Gamma,\Gamma,id)$ on $\mathfrak{sl}(n)$, 
one has
\begin{equation}
(\phi \o \phi)(\J_T) = 1 + 
\sum_{i < j} (q - q^{-1}) \frac{1}
{e^{\lambda_j - \lambda_i}-1} e_{ij} \o e_{ji},
\end{equation} 
where $\lambda_i$ denotes the $i$-th entry of $\lambda$.
\end{prop}

{\it Proof.}  This can be computed by 
using the product formula 
\begin{equation}
\J_T = \prod_{m = 1}^{\infty}
(\Ad e^{-m\lambda} \o 1)\Ad q^{(m-1) \Omega_\h}(\R_0)
\end{equation} 
and evaluating
in the vector 
representation, using the expression for $\R$ given in \cite{KhT}. \qs 

\begin{prop}  For the ``flip'' triple $(\Gamma,\Gamma,T)$ where 
$\g = \mathfrak{sl}(3)$,  one has $(\phi \o \phi)(\J_T) =J_T q^Z$ where
\begin{multline}
J_T(\mu) = 1 + (q-q^{-1})
\bigl[\frac{e^{-\mu}}{1-q^{-1}e^{-2\mu}} e_{12} \o e_{32}
+ \frac{qe^{-2\mu}}{1-qe^{-2\mu}} e_{12} \o e_{21} + 
\frac{e^{-\mu}}{1-qe^{-2\mu}}
e_{23} \o e_{21} \\ + \frac{q^{-1}e^{-2\mu}}{1-q^{-1}e^{-2\mu}} 
e_{23} \o e_{32} + \frac{-q^{-1} e^{-2\mu} + e^{-4\mu} + q^2e^{-4\mu} - 
qe^{-6\mu}}{(1-e^{-4\mu})(1-qe^{-2\mu})} e_{13} \o e_{31}\bigr],
\end{multline}
letting $\mu$ denote $\lambda_1 - \lambda_2 = \lambda_2 - \lambda_3$.
\end{prop}

{\it Proof.}  This can be seen by expanding 
\begin{equation}
\J_T q^{-Z} = \prod_{m=1}^{\infty}
\left(\Ad q^{[(T+...+T ^{m-1}) \o 1]\Omega_\h}(T^m\Ad e^{-m \lambda} \o 1) 
(\R_0)\right),
\end{equation} 
again using the formula for $\R$ found in \cite{KhT}. \qs

One may use these formulas to explicitly compute the image of the twisted
R-matrix in the vector representation.

\noindent
{\large Addresses:}

\noindent
{\bf Pavel Etingof:} Department of Mathematics, Room 2-165, MIT, 77 Massachusetts
Avenue,

\hskip 70 pt 
Cambridge, MA, 02139.

\noindent
{\bf Travis Schedler:} 059 Pforzheimer House Mail Center, Cambridge, MA, 02138.

\noindent
{\bf Olivier Schiffmann:} Department of Mathematics, MIT, 77 Massachusetts
Avenue, 

\hskip 101 pt
Cambridge, MA, 02139.

\noindent
{\large Email addresses:}

\noindent {\bf Pavel Etingof:} {\it etingof@math.harvard.edu}

\noindent {\bf Travis Schedler:} {\it schedler@fas.harvard.edu}

\noindent {\bf Olivier Schiffmann:} {\it schiffma@clipper.ens.fr}

\end{document}